\documentclass [12pt]  {article} 
\usepackage{amssymb}
\def \NN{{\mathbb{N}}}
\def \ZZ{{\mathbb{Z}}}

\def \CC{{\mathbb{C}}}
\def \PP{{\mathbb{P}}}
\def \DD{{\mathbb{D}}}

\begin{document}

\begin{center}
{\Large {\bf Value-sharing of meromorphic functions on a Riemann surface}}\\
\bigskip
{\sc Andreas Schweizer}\\
\bigskip
{\small {\rm National Center for Theoretical Sciences\\
Mathematics Division\\
National Tsing Hua University\\
Hsinchu 300, Taiwan\\
e-mail: schweizer@math.cts.nthu.edu.tw}}
\end{center}
\begin{abstract}
\noindent
We present some results on two meromorphic functions from $S$ to 
$\widehat{\CC}$ sharing a number of values where $S$ is a Riemann surface
of one of the following types: compact, compact minus finitely many points,
the unit disk, a torus, the complex plane.

\end{abstract}

\subsection*{Introduction}

We write $\widehat{\CC}$ (or sometimes $\PP^1$) for the Riemann sphere 
$\CC\cup\{\infty\}$. Let $R$ be a Riemann surface. Two meromorphic 
functions $f_1, f_2:R\to\widehat{\CC}$ are said to share the value 
$a\in\widehat{\CC}$ if for every $u\in R$ we have
$f_1(u)=a\Leftrightarrow f_2(u)=a$.
If moreover $f_1$ takes the value $a$ at each $u$ with the same
multiplicity as $f_2$, we say that $f_1$ and $f_2$ share the value
$a$ CM (counting multiplicities). If we don't know the multiplicities
or don't care, we say that $a$ is shared IM (ignoring multiplicities).
\par
In [Sa] Sauer proved among other results that if $S$ is a compact Riemann
surface of genus $g>0$ then two different non-constant meromorphic functions 
on $S$ cannot share more than $2+2\sqrt{g}$ values. If moreover one of the 
shared values is shared CM, the bound can be strengthened to
$\frac{1}{4}(9+\sqrt{32g+17})$.
In [Sch] these bounds have been slightly improved, and bounds in terms of 
other invariants of $S$ have also been given.
\par
In this paper we show that there do indeed exist compact Riemann surfaces 
that can carry two meromorphic functions with many shared values and that
one can even prescribe the shared values (Theorem 1). This question had 
been left open in [Sa] and [Sch]. The best examples one can find in the 
literature have $4$ shared values.
We also investigate how many values two functions that are meromorphic on 
a punctured compact Riemann surface can share (Section 3).
In the last section we use the concept of weighted sharing to refine
some known results.

\subsection*{1. Basic facts}

Recall that the gonality $d$ of a compact Riemann surface $S$ is defined 
to be the smallest integer $m$ such that $S$ can be realized as an 
$m$-sheeted branched covering of the Riemann sphere. Equivalently, $d$ is 
the smallest possible degree of a non-constant meromorphic function on $S$.
\\ \\
{\bf Theorem A.} [Sch] \it
Let $S$ be a compact Riemann surface of genus $g>0$ and gonality $d$.
Let $f_1, f_2: S\to\widehat{\CC}$ be two different non-constant 
meromorphic functions sharing $n$ values.
\begin{itemize}
\item[a)] Then 
$$n\le\min\{2+\sqrt{2g+2},\ 2d+1,\ 4+\frac{2(g-1)}{d}\}.$$
\item[b)] If moreover one of these $n$ values is shared CM, then we even have
$$n\le\min\{\frac{1}{2}(5+\sqrt{4g+5-2d}),\ 2d+1,\ 3+\frac{2(g-1)}{d}\}.$$
\item[c)] If $g=0$, i.e. $S=\PP^1$, the optimal bound in both cases is known 
to be $3$.
\end{itemize}
\rm
Actually, [Ba, Theorem 1] claims an inequality that is stronger than the
key lemmas in [Sa] and [Sch] and that would consequently imply stronger 
results than Theorem A. But trying to follow the proof in [Ba], I have 
not been able to obtain the claimed inequality.
\par
Besides the three articles already mentioned, the only other instance in the
literature dealing with value-sharing of meromorphic functions on a compact 
Riemann surface (other than the Riemann sphere) seems to be [A\&W]. It 
investigates functions sharing sets of values and is formulated in terms of 
functions on algebraic curves, which allows working over any algebraically 
closed field of characteristic $0$, not just over $\CC$. We point out that 
by the nature of their proofs the results in [Sa], [Sch] and some of the 
results we will obtain below (notably Theorem 1) also hold in such a general 
setting.
\par
The most famous result on value-sharing is of course Rolf Nevanlinna's 
theorem, that two meromorphic functions in the complex plane that share 
$5$ values must be equal. A true generalization of this result to Riemann
surfaces would, for example, be a statement about meromorphic functions on 
a punctured compact Riemann surface. I thank Jun-Muk Hwang who pointed 
this out to me and thus triggered my interest in this problem. Obviously,
this is more difficult than working on a compact Riemann surface, as there 
are many more meromorphic functions, and algebraic arguments cannot suffice. 
Luckily, the local result one would like to have has already been proved 
by R. Nevanlinna.
\\ \\
{\bf Theorem B.} [Ne] \it
Let $f_1(z)$ and $f_2(z)$ be meromorphic functions in a neighbourhood of
the point $z=\infty$, where they have an essential singularity. If for 
{\sc five} different values (finite or not) of $w$ the equalities
$$f_1(z)=w,\ \ \ \ f_2(z)=w,$$
outside some circle $|z|=r_0$ are satisfied for exactly the same values 
of $z$, then $f_1$ and $f_2$ are identical.
\rm
\\ \\
Of course one can transform this into a statement about meromorphic 
functions with an essential singularity in a punctured disk.
\par
Together with Theorem A c), which is folklore, Theorem B implies 
Nevanlinna's Five Value Theorem in the complex plane in the same way as 
the Big Picard Theorem implies the Little Picard Theorem. 
\\

\subsection*{2. Compact Riemann surfaces}

No example seems to be known of a compact Riemann surface $S$ and two
different nonconstant meromorphic functions $f_1$, $f_2:S\to\widehat{\CC}$
that share more than $4$ values. Moreover, [Sch] shows that ``most'' compact 
Riemann surfaces of a given genus do not allow more than $7$ shared values.
\par 
So it is legitimate to wonder whether in this case there is an absolute 
bound on the number of shared values, valid for all compact Riemann surfaces. 
The final remarks of [Sch] advocate this point of view, suggesting that 
a possible approach could be to prove the existence of a bound for the 
number of shared values of two different meromorphic functions on the open 
unit disk. The logical connection is immediate by pulling back the functions 
from the Riemann surface to its universal covering, which for $g\ge 2$ is 
the open unit disk.
\par
I am grateful to J\"org Winkelmann who showed to me that one can construct
different meromorphic functions on the open unit disk with any finite number 
of shared values. Although his examples do not come from compact Riemann 
surfaces, they convinced me that there probably is no uniform bound for all 
compact Riemann surfaces and that one should rather try to find examples of 
compact Riemann surfaces that allow many shared values. By an algebraic 
argument we will now construct such examples. In fact, one can even prescribe
the shared values.
\\ \\
{\bf Theorem 1.} \it
Let $a_1, a_2,\ldots,a_n\in\widehat{\CC}$ be $n$ different values with 
$n\ge 2$. 
\begin{itemize}
\item[a)] There exists a compact Riemann surface $S$ of genus $g\le n^2$
and two different non-constant meromorphic functions $f_1$ and $f_2$ from 
$S$ to $\widehat{\CC}$ that share the values $a_1,\ldots,a_n$. 
\item[b)] There exists a compact Riemann surface $S$ of genus 
$g\le 2n^2 -5n+3$ and two different non-constant meromorphic functions 
$f_1$ and $f_2$ from $S$ to $\widehat{\CC}$ that share the values 
$a_1,\ldots,a_{n-1}$ IM and the  value $a_n$ CM.
\end{itemize}
{\bf Proof.} \rm \\ 
a) After a M\"obius transformation we can assume that all the values
$a_i$ are finite and non-zero. We consider the polynomial
$$F(X,Y)=(X-Y)^{n+1}+Y\prod_{i=1}^n (Y-a_i)\prod_{i=1}^n (X-a_i).$$
As a polynomial in $X$ the highest coefficient is $1$, all other coefficients
are divisible by $Y$, and the absolute term is not divisible by $Y^2$.
(Here we are using that the $a_i$ are non-zero.) So it is an Eisenstein
polynomial in $X$ with respect to $Y$. By the Eisenstein Criterion
(see for example [Sti, Proposition III.1.14]) it is therefore irreducible
in $X$. Since the highest term is $X^{n+1}$, we also cannot factor out 
a polynomial that depends only on $Y$. So $F(X,Y)$ is irreducible.
\par
Let $S$ be the compact Riemann surface of the algebraic equation 
$$F(X,Y)=0.$$
Then the field of meromorphic functions on $S$ is $\CC(X,Y)$ where $X$ and
$Y$ are related by $F(X,Y)=0$. In particular, $X$ and $Y$ are functions of 
degree $n+1$ from $S$ to $\widehat{\CC}$. If $X$ takes the value $a_i$ at 
some point of $S$, we see from the equation that $Y$ cannot have a pole at 
that point and that it must take the same value. And vice versa. Thus $X$ 
and $Y$ share the value $a_i$.
\par
To estimate the genus of $S$ we apply the Castelnuovo Inequality. See 
[Sti, Theorem III.10.3] for an algebraic proof or [Ac, Theorem 3.5] for 
a proof that is more in the spirit of Riemann surfaces. The condition that 
the two maps do not factor over another Riemann surface corresponds to the
condition that $X$ and $Y$ generate the function field of $S$, which holds
by construction. From the Castelnuovo Inequality we get 
$$g\le(\deg(X)-1)(\deg(Y)-1)=n^2.$$
b) We can suppose that $a_n =\infty$ and that $a_1,\ldots, a_{n-1}$ are 
non-zero. This time we look at the Riemann surface corresponding to
$F(X,Y)=0$ where
$$F(X,Y)=(X-Y)^{2n-1}+Y\prod_{i=1}^{n-1}(Y-a_i)\prod_{i=1}^{n-1}(X-a_i).$$
As above we see that $F(X,Y)$ is irreducible and that the functions $X$ and 
$Y$ share the values $a_1,\ldots, a_{n-1}$. From the equation we also see
that if $X$ or $Y$ has a pole at some point, then the other function must 
also have a pole at that point with the same multiplicity.
\par 
Finally, $F(X,Y)=0$ is a plane curve of degree degree $d=2n-1$. Using the
formula for the genus of a plane curve ([Ac, p.5] or 
[Sti, Proposition III.10.5]), we get 
$$g\le\frac{(d-1)(d-2)}{2}=(n-1)(2n-3).$$
\hfill $\Box$
\\ \\
{\bf Remark 2.}
Together with the results from [Sa] (or Theorem A) this answers a question
asked in [Sa]. In both cases (all values shared IM, or one of the values
shared CM) the order of magnitude of the optimal bound in terms of the 
genus is $\sqrt{g}$. 
\par
There is still some room for further fine-tuning since the upper bounds 
from Theorem A grow asymptotically like $\sqrt{2g}$ resp. $\sqrt{g}$
whereas the examples from Theorem 1 grow asymptotically like $\sqrt{g}$
resp. $\sqrt{g/2}$.
\\ \\
{\bf Corollary 3.} \it
Let $n\in\NN$ and $a_1,a_2,\ldots,a_n\in\widehat{\CC}$. Then for every $k$ 
with $1\le k\le n$ there exist two different meromorphic functions $h_1$ and 
$h_2$ from the open unit disk $\DD=\{z\in\CC\ :\ |z|<1\}$ to $\widehat{\CC}$ 
that 
\begin{itemize}
\item omit the values $a_1,\ldots, a_{k-1}$;
\item take every other value in $\widehat{\CC}$ infinitely often;
\item share the values $a_k,\ldots, a_n$;
\item even share the value $a_n$ CM.
\end{itemize}
{\bf Proof.} \rm
Take a compact Riemann surface $S$ of genus $g\ge 2$ and two different 
functions $f_i:S\to\widehat{\CC}$ sharing $a_1,\ldots,a_n$ as in 
Theorem 1 b). Remove the inverse images of $a_1,\ldots,a_{k-1}$ from 
$S$, and restrict $f_i$ to functions $\widetilde{f_i}$ on the punctured 
Riemann surface $R$. Now we simply have to take 
$h_i=\widetilde{f_i}\circ\pi$ where $\pi:\DD\to R$ is the universal covering.
\hfill $\Box$
\\ \\
For a recent detailed study of value-sharing of meromorphic functions
in the unit disk see [Ti]. 
\\

\subsection*{3. Punctured Riemann surfaces}

Let $P_1$ be a point of $\PP^1$. By Nevanlinna's Five Value Theorem,
the maximum number of values that two different non-constant functions 
that are meromorphic on the punctured Riemann surface $\PP^1 -\{P_1 \}$ 
can share is $4$, realized for example by $e^z$ and $e^{-z}$ on $\CC$.
\par
The same bound holds for $\PP^1$ minus two points. After a M\"obius 
transformation we can assume that $f_1$ and $f_2$ are meromorphic on
$\CC -\{0\}$. Then $f_1(e^z)$ and $f_2(e^z)$ are meromorphic on $\CC$
and share the same values; so their number is bounded by $4$.
\\ \\
{\bf Examples 4.}
\begin{itemize}
\item[a)]
Let $\zeta_r =e^{\frac{2\pi i}{r}}$ be a primitive $r$-th root of unity
with $r\ge 2$. Then the functions 
$$f_1(z)=z\ \ \ \hbox{\rm and}\ \ \ f_2(z)=\zeta_r z$$ 
are obviously meromorphic on 
$$\PP^1 -\{\zeta_r^k\ :\ k=0,\ldots,r-1\}$$
and share the $r+2$ values
$$\infty,\ 0 \ \ \hbox{\rm and}\ \ \zeta_r^k,\ \ (k=0,\ldots,r-1).$$
\item[b)]
For any fixed choice of three points $P_1, P_2, P_3$ on $\PP^1$ we can 
always construct two meromorphic functions on $\PP^1 -\{P_1,P_2,P_3\}$
that share $5$ values. Let $T(z)$ be the M\"obius transformation that maps
$P_1, P_2, P_3$ to $1, \zeta_3, \zeta_3^2$. Then $T(z)$ and $\zeta_3 T(z)$
share $\infty, 0, 1, \zeta_3, \zeta_3^2$, actually all CM.
\item[c)]
On any fixed Riemann surface $\PP^1 -\{P_1,\ldots,P_r\}$ with $r\ge 4$
we can realize at least $6$ shared values. Let $T_1(z)$ be the M\"obius
transformation that maps $P_1, P_2, P_3$ to $\infty, 0, 1$ and let 
$c=T_1(P_4)$. Fix a square-root $\sqrt{c}$ and set 
$T_2(z)=\frac{z+\sqrt{c}}{z-\sqrt{c}}$. Then $T_2$ maps $\infty, 0, 1, c$
to $1, -1, a, -a$. So the functions $T_2(T_1(z))$ and $-T_2(T_1(z))$
from $\PP^1 -\{P_1,\ldots,P_r\}$ to $\widehat{\CC}$ share
$\infty, 0, 1, -1, a, -a$, actually all CM.
\item[d)]
If $g\ge 2$ and $R_1,\ldots,R_{2g+2}\in\CC$ are different, the compact 
Riemann surface $S$ corresponding to 
$$Y^2 =\prod_{i=1}^{2g+2}(X-R_i)$$
is hyperelliptic of genus $g$. Moreover, the hyperelliptic map
$\kappa:S\to\PP^1$, corresponding to $(X,Y)\mapsto X$, is ramified 
exactly above the points $R_i$. If $\{ R_1,\ldots,R_{2g+2}\}$ contains
all $r$-th roots of unity and $P_i =\kappa^{-1}(\zeta_r^i)$ for
$i=1,2,\ldots,r$, then $\kappa$ and $\zeta_r\kappa$ share the $r+2$
values $\infty$, $0$, $1$, $\zeta_r,\ldots,\zeta_r^{r-1}$ on
$S-\{P_1,\ldots,P_r\}$.
\end{itemize}
In Examples 4 a), b) and d) we actually have constructed $r$ different 
functions that all share the same $r+2$ values CM. Simple as these 
constructions may be, under certain conditions they will turn out to 
be essentially the only ones that realize the maximal possible number 
of shared values.
\\ \\
Now we combine Theorem B and algebraic arguments to derive upper bounds.
Recall that if $S$ is a compact Riemann surface and $P_1,\ldots,P_r$ are $r$ 
different points on $S$ with $0\le r<\infty$, then the Euler characteristic 
of the Riemann surface $R=S-\{P_1,\ldots,P_r\}$ is defined to be
$$\chi(R)=2-2g(S)-r.$$
If $\chi(R)<0$ then $R$ is called hyperbolic. For $r>0$ this only excludes 
the two cases $g=0$, $r\le 2$, which we have discussed at the beginning of
this section.
\\ \\
{\bf Theorem 5.} \it
Let $S$ be a compact Riemann surface of genus $g$ and gonality $d$.
Let $R=S-\{P_1,\ldots,P_r\}$ where $P_1,\ldots,P_r$ are $r$ different
points on $S$.
\par
If $R$ is hyperbolic, then two different non-constant meromorphic functions 
on $R$ can share at most 
$$4+\frac{2g-2+r}{d}=4-\frac{\chi(R)}{d}$$
values.
\par
Moreover, if $4-\frac{\chi(R)}{d}$ is an integer and if $f_1$ and $f_2$ 
realize this bound, then $f_1$ and $f_2$ must both be meromorphic on the 
compact Riemann surface $S$ and $\deg(f_1)=\deg(f_2)=d$.
\\ \\
{\bf Proof.} \rm
First let us assume that one of the functions, say $f_1$ has an essential 
singularity at one of the points $P_i$. By the Big Picard Theorem, in
every neighbourhood of $P_i$ the function $f_1$ can omit at most two
of the shared values. Every other shared value is a limit of $f_2(z)$
when $z$ approaches $P_i$ along the corresponding inverse images. So
if there are $4$ or more shared values, then $f_2$ also has an essential 
singularity at $P_i$. By Theorem B in this case $f_1$ and $f_2$ cannot 
share more than $4$ values.
\par
Now suppose that none of the points $P_i$ is an essential singularity.
Then $f_1$ and $f_2$ extend to meromorphic functions on the compact
Riemann surface $S$. Let $d_i =\deg(f_i)$. Without loss of generality 
we can assume $d_1\le d_2$. After a M\"obius transformation we can also 
assume that all $n$ shared values $a_1,\ldots,a_n$ lie in $\CC$. Let 
$M$ consist of all $u\in S$ with $f_2(u)\in\{a_1,\ldots,a_n\}$. Then 
$$r_2(M):=\sum_{u\in M}(\hbox{\rm mult}_{f_2}(u)-1)$$
measures the ramification of the covering $f_2:S\to\widehat{\CC}$ above 
these values. Applying the Hurwitz formula we get 
$$r_2(M)\le r_2(S)=2g(S)-2-d_2(2g(\widehat{\CC})-2)=2g-2+2d_2.$$
Also, every element of $M\cap R$ is a zero of $f_1 -f_2$, which is 
a function of degree $\le 2d_2$. Together we obtain
$$nd_2=|M\cap R|+|M\cap\{P_1,\ldots,P_r\}|+r_2(M)\le 2d_2 +r+2g-2+2d_2,$$
and after division
$$n\le 4+\frac{2g-2+r}{d_2}=4-\frac{\chi(R)}{d_2}.$$
Since (by definition) $d\le d_1\le d_2$, this establishes the bound, and
it also shows that reaching the bound is only possible if $d_2 = d$.
\hfill $\Box$
\\ \\
{\bf Corollary 6.} \it
Let $P_1,\ldots, P_r$ be $r$ different points on $\PP^1$ with $r\ge 3$.
Then two different non-constant functions that are meromorphic on 
$\PP^1 -\{P_1,\ldots,P_r\}$ cannot share more than $r+2$ values.
\par
Moreover, this bound is optimal (at least for a suitable choice of
$P_1,\ldots,P_r$). If $f_1$ and $f_2$ attain this bound, they both must 
be functions of degree $1$ on $\PP^1$, that is, they must be fractional 
linear transformations.
\\ \\
{\bf Proof.} \rm
Specialize Theorem 5. The bound is sharp by Example 4 a).
\hfill $\Box$
\\ \\
Note that if $r>4$ we do not claim that the bound in Corollary 6 is sharp 
for every choice of $P_1,\ldots,P_r$.
\\ \\
{\bf Example 7.}
If the $5$ points $P_1,\ldots,P_5$ from $\PP^1$ are not all lying on one
circle or on one straight line, then two different non-constant meromorphic
functions on $R=\PP^1 -\{P_1,\ldots,P_5\}$ cannot share more than $6$ values.
\par
Indeed, if $f_1$ and $f_2$ share $7$ values, then by Corollary 6 they must
be functions of degree $1$ on $\PP^1$. Moreover, two of the shared values 
must obviously be taken at points of $R$. Applying M\"obius transformations
to the values and to the arguments we can assume that these two shared
values are $\infty$ (taken at the point $\infty$) and $0$ (taken at the
point $0$). Hence $f_i(z)=c_i z$ with $c_i\in\CC^*$. Without loss of 
generality we can assume $f_1(z)=z$. Then $f_2$ must permute the $5$ 
punctures. This is only possible if $f_2$ is a rotation around $0$ and
if the $5$ punctures are lying on a circle around $0$. Since M\"obius
transformations respect circles on $\PP^1$, the $5$ original points
$P_1,\ldots,P_5$ must lie on such a circle.
\\ \\
For $r=0$ Theorem 5 gives one of the bounds from Theorem A. In view of the
other bounds in Theorem A one might think that in the case where $d$ is small 
and $g$ is large it should be possible to get stronger bounds than Theorem 5.
But even under this condition the bound in Theorem 5 is often sharp.
\\ \\
{\bf Proposition 8.} \it
For every compact Riemann surface $S$ there are infinitely many $r\in\NN$
for which the bound in Theorem 5 is sharp (provided the points 
$P_1,\ldots,P_r$ are suitably chosen), and the values are even shared CM.
\\ \\
{\bf Proof.} \rm
Fix a covering $\pi:S\to\PP^1$ of degree $d$. We can assume that all
ramified values $R_1,\ldots,R_m$ of $\pi$ are in $\CC^*$.
\par 
If $Q_1,\ldots,Q_s$ are $s$ different points in $\CC^*$, containing all 
$R_i$, then by the Hurwitz formula there are exactly $r=ds-(2d+2g-2)$
points $P_i$ of $S$ lying above $Q_1,\ldots,Q_s$. 
\par
If moreover the set $\{Q_1,\ldots,Q_s\}$ is closed under $Q_i\mapsto -Q_i$,
then the functions $\pi$ and $-\pi$ from $S-\{P_1,\ldots,P_r\}$ to 
$\widehat{\CC}$ share the $2+s$ values $\infty$, $0$, $Q_1,\ldots,Q_s$
CM and $2+s=4+\frac{2g-2+r}{d}$.
\hfill $\Box$
\\ \\
However, for $r$ in a certain range one can indeed improve on Theorem 5.
\\ \\
{\bf Theorem 9.} \it
Let $S$ be a compact Riemann surface of genus $g$ and gonality $d$.
Let $R=S-\{P_1,\ldots,P_r\}$ where $P_1,\ldots,P_r$ are $r$ different
points on $S$.
\par
If $r\ge 2d$, then two different non-constant meromorphic functions 
on $R$ can share at most $r+2$ values.
\par
Moreover, if $f_1$ and $f_2$ realize this bound and $r>2d$, then $f_1$ 
and $f_2$ must both be meromorphic on $S$ with $\deg(f_1)=\deg(f_2)=d$,
and the coverings $f_i :S\to\PP^1$ are totally ramified at all points
$P_1,\ldots,P_r$.
\\ \\
{\bf Proof.} \rm
As explained at the beginning of the proof of Theorem 5, we can assume
that $f_1$ and $f_2$ are meromorphic on $S$. Without loss of generality
$\deg(f_1)=d_1 \le d_2 =\deg(f_2)$.
\par
If $(d-1)(d_2 -1)\ge g$, i.e. $d_2\ge\frac{g-1+d}{d-1}$, from the last
inequality in the proof of Theorem 5 we get
$$n\le 4+\frac{2g-2+r}{d_2}\le 4+\frac{(2g-2+r)(d-1)}{g-1+d}
=4+2(d-1)+\frac{(r-2d)(d-1)}{g-1+d},$$
which is smaller than $2+2d+r-2d=r+2$ if $r>2d$ and equal to $r+2$ if
$r=2d$.
\par
If $(d-1)(d_2 -1)<g$, we proceed by induction on $d$, using 
Corollary 6 as induction basis.
\par
If $d>1$ let $F$ be the function field of $S$. Fix a rational subfield 
$R$ of index $d$ in $F$. Let $M$ be the compositum of $R$ and 
$\CC(f_2)$ and let $c$ be the index of $M$ in $F$. 
By Castelnuovo's inequality we have
$$\textstyle{g(M)\le(\frac{d}{c}-1)(\frac{d_2}{c}-1).}$$ 
Now $c=1$ would be equivalent to $M=F$ and hence contradict the 
Castelnuovo inequality. So $M$ is a proper subfield of $F$.
\par
If the compositum of $M$ and $\CC(f_1)$ were $F$, then again
by Castelnuovo's inequality we would obtain 
$g\le c g(M)+(c-1)(d_1 -1)$, and hence 
$${\textstyle (d-1)(d_2-1)<g\le c(\frac{d}{c}-1)(\frac
{d_2}{c}-1)+(c-1)(d_2-1).}$$
Subtracting $(c-1)(d_2 -1)$ we get the contradiction 
$(d-c)(d_2 -1)<(d-c)(\frac{d_2}{c}-1)$.
We conclude that $R$, $\CC(f_1)$ and $\CC(f_2)$ are contained 
in a proper subfield $L$ of $F$. Let $\delta=[F:L]$. Then 
$$f_i =\widetilde{f_i}\circ\kappa$$ 
where $\widetilde{f_1}$ and $\widetilde{f_2}$ are meromorphic of degrees
$\frac{d_1}{\delta}$ resp. $\frac{d_2}{\delta}$ on the compact Riemann
surface $\widetilde{S}$ corresponding to $L$ and $\kappa$ is the 
covering map from $S$ to $\widetilde{S}$.
\par
Now the image of $\{ P_1,\ldots,P_r\}$ under $\kappa$ is a subset of
$\widetilde{S}$ of cardinality $\widetilde{r}$ where obviously
$$\frac{r}{\delta}\le\widetilde{r}\le r.$$ 
Moreover, $\widetilde{r}\ge 2\frac{d}{\delta}$ and $\frac{d}{\delta}$ 
is the gonality of $\widetilde{S}$.
By induction $\widetilde{f_1}$ and $\widetilde{f_2}$ can share at most
$\widetilde{r}+2$ values, and if they share $\widetilde{r}+2$ values
then $\deg(\widetilde{f_i})=\frac{d}{\delta}$ and $\widetilde{f_1}$ and
$\widetilde{f_2}$ are totally ramified at all $\kappa(P_i)$.
\par
This proves the first statement of the theorem. It also shows that
having $r+2$ shared values is only possible if $\widetilde{r}=r$,
i.e. if $\kappa$ is totally ramified at $P_1,\ldots,P_r$.
\hfill $\Box$
\\ \\
{\bf Remarks 10.} 
\begin{itemize}
\item[a)] 
In the range $2d\le r<\frac{2g}{d-1}+2$ the bound in Theorem 9 is 
better than Theorem 5 if $d>1$. Of course, this range might be empty
if $d$ is sufficiently big with respect to $g$.
\item[b)]
By Example 4 d), for each pair $(g,r)$ with $4\le r\le 2g+2$ there exists 
a hyperelliptic Riemann surface $S$, points $P_1,\ldots,P_r$ on $S$, and
meromorphic $f_1,f_2$ on $S-\{P_1,\ldots,P_r\}$ sharing $r+2$ values.
\par
On the other hand, if $S$ is the hyperelliptic surface corresponding to
$$Y^2=(X-1)(X-2)(X-3)(X-1-i)(X-2-i)(X-3-i),$$
then by Theorem 9 and Example 7 for every choice of five points 
$P_1,\ldots,P_5$ on $S$ we cannot get more than $6$ shared values.
\item[c)]
More generally, if $d\ge 3$ and $g=\frac{1}{2}(d-1)(md-2)$ with
$m\ge 2$, then for every $r$ with $2d\le r\le md=\frac{2g}{d-1}+2$
there exists a compact Riemann surface $S$ of genus $g$ and gonality 
$d$, points $P_1,\ldots,P_r$ on $S$, and meromorphic $f_1$ and $f_2$
on $S-\{P_1,\ldots,P_r\}$ that share $r+2$ values.
\par
Explicitly, let $S$ be the Riemann surface of the function field
$F=\CC(X,Y)$ with $Y^d =f(X)$ where $f(X)\in\CC[X]$ is square-free,
of degree $md$ and divisible by $X^r -1$. Using the Hurwitz formula
one can calculate the genus of every intermediate field between 
$\CC(X)$ and $F$. Then the Castelnuovo inequality shows that $F$ 
indeed has gonality $d$. Let $\pi$ be the covering map from $S$ to
$\PP^1$ corresponding to the extension $F/\CC(X)$. Let
$P_i =\pi^{-1}(\zeta_r^i)$ for $i=1,2,\ldots,r$. Then $\pi$ and 
$\zeta_r\pi$ share the values 
$\infty$, $0$, $1$, $\zeta_r,\ldots,\zeta_r^{r-1}$ 
on $S-\{P_1,\ldots,P_r\}$.
\end{itemize}
We finish this section with another example of which one can easily 
construct many explicit instances.
\\ \\
{\bf Example 11.} 
Let $F$ be the compositum of two quadratic extensions of $\CC(z)$. Then
$F$ is a Galois extension of $\CC(z)$ with Galois group 
$\ZZ/2\ZZ\oplus\ZZ/2\ZZ$. Assume that $F$ has genus $g\ge 10$ and that
each of the three intermediate quadratic fields has genus at least $2$.
Then the Castelnuovo inequality implies that the gonality of $F$ is $4$
and that $\CC(z)$ is the only rational subfield over which $F$ has degree 
$4$. Moreover, no place is totally ramified in $F/\CC(z)$. 
\par
Let $S$ be the compact Riemann surface of $F$. Then Theorem 9 implies 
that for any choice of $r>8$ different points $P_1,\ldots,P_r$ on $S$ 
two meromorphic functions on $S-\{P_1,\ldots,P_r\}$ cannot share more 
than $r+1$ values.
\\

\subsection*{4. Weighted sharing}

In order to refine the results we recall the concept of weighted sharing 
as introduced by Lahiri in [La]. A shared value $a$ is shared with weight
$m\in\NN_0\cup\{\infty\}$ if for all inverse images $u$ of $a$ we have
$$\hbox{\rm mult}_{f_1}(u)=\mu\le m \Leftrightarrow
\hbox{\rm mult}_{f_2}(u)=\mu\le m$$
and
$$\hbox{\rm mult}_{f_1}(u)>m \Leftrightarrow\hbox{\rm mult}_{f_2}(u)>m.$$
In particular, $f_1$ and $f_2$ sharing the value $a$ with weight one means 
that the simple $a$-points of $f_1$ are exactly the simple $a$-points of
$f_2$ and the multiple $a$-points of $f_1$ are exactly the multiple 
$a$-points of $f_2$, where in the latter case the multiplicities are not
necessarily the same.
\par
Obviously, sharing with weight $0$ simply means sharing IM, and
sharing with weight $\infty$ is the same as sharing CM.
\par
In the sequel we write $(m_1 ,\ldots ,m_n)$ to indicate that the value
$a_i$ is shared with weight $m_i$.
\par
Since a meromorphic function on a compact Riemann surface $S$ is determined
up to  a multiplicative constant by its divisor, we cannot have two different 
non-constant meromorphic functions on $S$ sharing $3$ values with weights
$(\infty,\infty,0)$. (Apply a M\"obius transformation to move the two
CM-shared values to $0$ and $\infty$.) But for every $m\in\NN$ sharing 
with weights $(\infty,m,0)$ is possible on every $S$.
\\ \\
{\bf Example 12.}
From [Pi] we take the example of the functions 
$$f_1(z)=\frac{-4z^3}{(z-1)^3(z+1)}\ \ \hbox{\rm and }\ \ 
f_2(z)=\frac{-4z}{(z-1)(z+1)^3}$$
from $\widehat{\CC}$ to $\widehat{\CC}$. They share the value $1$ CM
(taken with multiplicity $1$ at the zeroes of $(z^2 +1)(z^2 +2z-1))$, 
and the values $0$ and $\infty$, both IM. Since the value $0$ is taken 
exactly at the points $0$ and $\infty$, the functions $h_i(z)=f_i(z^{m+1})$ 
take this value at these two points with multiplicity at least $m+1$. 
So, somewhat trivially, they share the value $0$ with weight $m$.
\par
Finally, if $S$ is any compact Riemann surface, there is a covering
$\pi:S\to\PP^1$. Then $h_1 \circ\pi$ and $h_2 \circ\pi$ are two meromorphic
functions on $S$ that share the values $1$, $0$ and $\infty$ with 
respective weights $(\infty,m,0)$. Applying a suitable M\"obius 
transformation, we can get any $3$ values shared with these weights.
\\ \\
{\bf Remark 13.}
If $a_1, a_2,\ldots,a_n\in\widehat{\CC}$ and $m\in\NN$, with the same trick 
we can even construct a compact Riemann surface $R$ and two functions $h_1$,
$h_2$ from $R$ to $\widehat{\CC}$ that share the values $a_1,\ldots,a_n$ with
respective weights $(\infty, m, m, \ldots, m)$. 
We only have to take $S$ and $f_1$, $f_2$ as in Theorem 1 b), take a finite 
covering $\pi:R\to S$ of compact Riemann surfaces that is totally ramified 
above all inverse images of the shared values with ramification index at 
least $m+1$, and set $h_i=f_i\circ\pi$.
\par
The same weighted sharing can then of course also be obtained in Corollary 3.
\\ \\
Admittedly, in the preceeding examples the weight $m$ of the sharing does
not really tell much, as we simply have artificially increased the
multiplicities. However, if one restricts the nature of the underlying
Riemann surface, one can get non-trivial information, as we will show now 
by improving a result from [Sch]. 
\par
The following lemma can be easily obtained as a special case of known
results on two meromorphic functions in the complex plane that share
$4$ values (cf. [Gu1], [Gu2], [Mu], [Y\&Y]). We prefer to give a direct 
algebraic proof. Besides being conceptually much simpler, it has the 
advantage to be valid for rational functions on an elliptic curve over 
any algebraically closed field of characteristic $0$.
\\ \\
{\bf Lemma 14.} \it
Let $S$ be a compact Riemann surface of genus $1$ and $f_1$, $f_2$ two 
non-constant meromorphic functions from $S$ to $\widehat{\CC}$ that share 
the $4$ values $a_1,a_2,a_3,a_4 \in\widehat{\CC}$. Let 
$M=f_i^{-1}\{a_1,a_2,a_3,a_4\}\subseteq S$. Then
\begin{itemize}
\item[a)] $\deg(f_1)=\deg(f_2)$; 
\item[b)] $|M|=2\deg(f_i)$;
\item[c)] the map $f_i :S\to\widehat{\CC}$ is unramified outside $M$;
\item[d)] $f_1(u)\neq f_2(u)$ for every $u\not\in M$;
\item[e)] at each point $u\in M$ at least one of the two functions takes
the shared value $a_j =f_i(u)$ with multiplicity $1$.
\end{itemize}
{\bf Proof.} \rm
Let $d_i =\deg(f_i)$. Without loss of generality we can assume $d_1\le d_2$.
Let 
$$r_2(M):=\sum_{u\in M}(\hbox{\rm mult}_{f_2}(u)-1).$$
Applying the Hurwitz formula to the covering $f_2:S\to\widehat{\CC}$ we get 
$$r_2(M)\le r_2(S)=2g(S)-2-d_2(2g(\widehat{\CC})-2)=2d_2.$$
After a M\"obius transformation we can assume $a_1,\ldots,a_4\in\CC$.
Since every $u\in M$ is a zero of $f_1 -f_2$, we have
$$d_1 +d_2\ge\deg(f_1 -f_2)\ge |M|=4d_2-r_2(M)\ge 2d_2.$$
This shows $d_1 =d_2$ and $|M|=2d_2$ and also $r_2(M)=2d_2$, which means 
that $f_2$ is unramified outside $M$. Since $d_1 =d_2$ we can interchange
$f_1$ and $f_2$, so $f_1$ is also unramified outside $M$.
\par
Moreover, we see $\deg(f_1 -f_2)=2d_2$. This shows that $f_1$ and $f_2$ 
have no common poles. Finally, since $f_1 -f_2$ vanishes at the $2d_2$ 
different points in $M$, it cannot vanish outside $M$ (so claim d) holds) 
and it cannot have a multiple zero, which implies statement e).
\hfill $\Box$
\\ \\
{\bf Corollary 15.} \it
Let $S$ be a compact Riemann surface of genus $1$ and $f_1$, $f_2$ two 
non-constant meromorphic functions from $S$ to $\widehat{\CC}$ that share 
$4$ values with respective weights $(1,0,0,0)$. Then $f_1 =f_2$.
\rm
\\ \\
{\bf Proof.} \rm
By Lemma 14 e) the value that is shared with weight one is actually shared
CM. By part b) of Theorem A this implies $f_1= f_2$.
\hfill $\Box$
\\ \\
{\bf Corollary 16.} \it
Let $f_1$ and $f_2$ be two non-constant elliptic functions on the complex 
plane (not necessarily with commensurable period lattices). If $f_1$ and 
$f_2$ share $4$ values, of which one is shared with weight one, then the
functions $f_1$ and $f_2$ are equal.
\rm
\\ \\
{\bf Proof.} \rm
Let $\Lambda_i$ be the period lattice of $f_i$. Translating the variable 
$z$ we can assume that $f_1(0)$ is one of the shared values. Then $f_1$ 
(and also $f_2$) takes this value at all points of the $\ZZ$-module 
$\Lambda_1 +\Lambda_2$. Since $f_1$ is not constant, $\Lambda_1 +\Lambda_2$ 
must be discrete and hence a rank $2$ lattice. Thus $\Lambda_1$ and 
$\Lambda_2$ are commensurable. Let $\Lambda$ be the rank $2$ lattice
$\Lambda_1\cap\Lambda_2$. Then we can consider $f_1$ and $f_2$ as meromorphic 
functions on the torus $\CC/\Lambda$ and apply the previous corollary.
\hfill $\Box$
\\ \\
This result, presumably well known to specialists, has inspired the 
following modification of a famous problem.
Let $f_1$ and $f_2$ be two non-constant meromorphic functions in the 
complex plane sharing $4$ values. Gundersen [Gu1] has shown that if $2$
of these $4$ values are shared CM then all $4$ values must be shared CM.
He also asked whether one CM-shared value would already be enough for the
same conclusion. Although there are positive answers under different 
additional conditions (see for example [Mu], [Gu2], [Y\&Y, Chapter 4],
[Hu]), this is still an open problem.
\par
But the notion of weighted sharing opens up infinitely many more 
possibilities between weights $(\infty, 0,0,0)$ and $(\infty,\infty,0,0)$. 
So one might ask: Does sharing with weights $(\infty, m,0,0)$ for some 
$m\in\NN$ imply weights $(\infty,\infty,\infty,\infty)$?
It turns out that a much weaker condition, namely that the weights of 
sharing are $(1,1,0,0)$, already suffices.
\\ \\
{\bf Theorem 17.} \it
Let $f_1$ and $f_2$ be two non-constant meromorphic functions in the complex 
plane sharing $4$ values. If $2$ of these values are shared with weight one, 
then all $4$ values are shared CM.
\\ \\
{\bf Proof.} \rm
If $f_1$ and $f_2$ share $4$ values, then by [Mu, Lemma 1] for each 
shared value $a_j$ the counting function of the points that are multiple 
$a_j$-points for both functions is $S(r,f_i)$. Thus a value that is shared 
with weight one is shared ``CM'' in the sense of [Mu], that is, the counting 
function of the points where the value is taken with different multiplicities 
is $S(r,f_i)$. But by [Gu2, Theorem C*], two ``CM''-shared values imply that 
all four values are shared CM.
\hfill $\Box$
\\ \\
Modifying the condition that one value is shared CM in the other direction 
by relaxing it, we obtain the following problem.
\\ \\
{\bf Question:} 
Given $m\in\NN$. If two non-constant meromorphic functions in the complex 
plane share $4$ values with respective weights $(m,0,0,0)$, does this imply
that all four values are shared CM?
\\ \\
This is presumably a very difficult question. A positive answer would settle
the famous problem on $1$ CM plus $3$ IM shared values. And a counterexample
would be a new example of two meromorphic functions in the complex plane 
sharing $4$ values not all of which are shared CM. There are essentially 
only $3$ known examples of such functions, namely the ones described in 
[Gu1], [Ste] and [Re].
\\ \\
{\bf Acknowledgements.} I thank Jun-Muk Hwang for the initial impulse,
J\"org Winkelmann for the decisive hint, and Andreas Sauer for several
helpful comments and discussions.

\subsection*{\hspace*{10.5em} References}
\begin{itemize}
\item[{[Ac]}] \sc R. Accola: \it Topics in the Theory of Riemann Surfaces,
\rm Springer Lecture Notes in Mathematics 1595, 
Berlin-Heidelberg-New York, 1994
\item[{[A\&W]}] \sc T.T.H. An and J.T.-Y. Wang: \rm Unique range 
sets and uniqueness polynomials for algebraic curves, \it 
Trans. Amer. Math. Soc. \bf 359 \rm (2007), no. 3, 937-964
\item[{[Ba]}] \sc E. Ballico: \rm Meromorphic functions on compact 
Riemann surfaces and value sharing, \it Complex Variables 
\bf 50, No. 15 \rm (2005), 1163-1164
\item[{[Gu1]}] \sc G. Gundersen: \rm Meromorphic functions that 
share four values, \it Trans. Amer. Math. Soc. \bf 277 \rm (1983), 
545-567\\
Correction: \it Trans. Amer. Math. Soc. \bf 304 \rm (1987), 847-850
\item[{[Gu2]}] \sc G. Gundersen: \rm Meromorphic functions that share 
three values IM and a fourth value CM, 
\it Complex Variables \bf 20 \rm (1992), 99-106
\item[{[Hu]}] \sc B. Huang: \rm On Gundersen's question for the unicity
of meromorphic functions, \it Anal. Theory Appl. \bf 21 \rm (2005), 
235-241
\item[{[La]}] \sc I. Lahiri: \rm Weighted sharing and uniqueness of
meromorphic functions, \it Nagoya Math. J. \bf 161 \rm (2001), 193-206
\item[{[Mu]}] \sc E. Mues: \rm Meromorphic Functions Sharing Four Values,
\it Complex Variables \bf 12 \rm (1989), 169-179
\item[{[Ne]}] \sc R. Nevanlinna: \rm Un th\'eor\`eme d'unicit\'e relatif
aux fonctions uniformes dans le voisinage d'un point singulier essentiel,
\it C. R. de l'Acad. des Sc. \bf 181 \rm (1925), 92-94
\item[{[Pi]}] \sc A. Pizer: \rm A problem on rational functions, 
\it Amer. Math. Monthly \bf 80 \rm (1973), 552-553
\item[{[Re]}] \sc M. Reinders: \rm A New Example of Meromorphic 
Functions Sharing Four Values and a Uniqueness Theorem, \it 
Complex Variables \bf 18 \rm (1992), 213-221
\item[{[Sa]}] \sc A. Sauer: \rm Uniqueness theorems for holomorphic 
functions on compact Riemann surfaces, \it New Zealand J. Math.
\bf 30 \rm (2001), 177-181
\item[{[Sch]}] \sc A. Schweizer: \rm Shared values of meromorphic 
functions on compact Riemann surfaces, \it Arch. Math. (Basel) 
\bf 84 \rm (2005), 71-78
\item[{[Ste]}] \sc N. Steinmetz: \rm A uniqueness theorem for three 
meromorphic functions, \it Ann. Acad. Sci. Fenn. Ser. A I Math.
\bf 13 \rm (1988), 93-110
\item[{[Sti]}] \sc H. Stichtenoth: \it Algebraic Function Fields 
and Codes, \rm Springer Universitext, Berlin-Heidelberg-New York,
1993
\item[{[Ti]}] \sc F. Titzhoff: \rm Slowly growing functions sharing 
values, \it Fiz. Mat. Fak. Moksl. Sem. Darb. \bf 8 \rm (2005), 143-164
\\
electronic access via \tt http://siauliaims.su.lt/article.al?id=29 \rm
\item[{[Y\&Y]}] \sc C.-C. Yang and H.-X. Yi: \it Uniqueness theory
of meromorphic functions, \rm Kluwer Academic Publishers Group,
Dordrecht, 2003
\end{itemize}

\end{document}